\documentclass[11pt]{amsart}

\theoremstyle{remark}

\theoremstyle{definition}

\numberwithin{equation}{section}

\begin{document}

\title{Dimension of  spaces of polynomials on  abelian topological  semigroups }
\author{Bolis Basit}

       \address {School of Math. Sciences, P.O. Box  28M,
       Monash University, Victoria 3800, Australia}
        \email{ bolis.basit@sci.monash.edu.au and alan.pryde@sci.monash.edu.au}
\author{A. J. Pryde}

        \begin{abstract}
{In this paper we study (continuous) polynomials   $p: J\to X$,
where $J$ is an abelian  topological   semigroup and $X$ is a
topological vector space.  If $J$ is a subsemigroup with non-empty
interior of a locally compact abelian group $G$ and $G=J-
 J$, then every
polynomial $p$ on $J$ extends uniquely to a polynomial on $ G$. It
is of particular interest to know when the spaces  $P^n (J,X)$ of
polynomials of order at most  $n$ are finite dimensional. For
example we show that for some semigroups the subspace $P^n_{R}
(J,\mathbf{C})$ of  Riss polynomials (those generated by a finite
number of homomorphisms $\alpha: J\to \mathbf{R}$)  is properly
contained in $P^n (G,\mathbf{C})$. However, if $P^1
(J,\mathbf{C})$ is finite dimensional then $P^n_{R}
(J,\mathbf{C})= P^n (J,\mathbf{C})$. Finally we exhibit a large
family of groups for which  $P^n (G,\mathbf{C})$ is  finite
dimensional.
\smallskip


 }
 \end{abstract}
\maketitle

{\center\textbf{0. Introduction and notation}\endcenter}

Throughout this paper $G$ will denote a locally compact abelian
group with dual (character) group $\widehat{G}$.  In developing a
Laplace transform for such groups, Mackey\cite{MWG} described the
relationship between the continuous homomorphisms $\alpha: G\to$
$\mathbf{R}$ and the continuous one parameter subgroups $\beta:
\mathbf{R} \to \widehat{G}$. The latter are precisely the duals
$\hat{\alpha}$ of the former, defined by $\hat{\alpha}(u)(t)=
e^{i\, u \alpha(t)}$ for $u\in\mathbf{R}$ and $t\in G$. Hence
there is a non-trivial (unbounded) homomorphism $\alpha: G\to
\mathbf{R}$ if and only if there is a non-trivial one parameter
subgroup $\beta: \mathbf{R} \to \widehat{G}$.

Riss \cite{RJ} introduced  a space, say $P_{R} (G)$, of
complex-valued polynomials on $G$. It is the  unital subalgebra
generated by the homomorphisms $\alpha: G\to \mathbf{R}$ in the
algebra $C(G)$ of continuous  complex-valued functions (in this
paper  homomorphisms are always continuous). Subsequently
Domar\cite{D} defined a
 different space of polynomials ($P(G)$ below) and mentioned that
 the connection between
 the two classes is not obvious and  a study of this problem would involve extensive structural
 considerations. See [4, end of Chapter II]. To the best of our knowledge
 this problem has not yet been resolved and this note is devoted in part to clarifying this connection.

Throughout $X$ will denote a topological vector  space over
$\mathbf {C}$  and $J$ an abelian topological semigroup. Replacing
$G$ by $J$ in the definition of Riss polynomials we obtain the
space $P_{R} (J)$.
 Following \cite{D} (see also \cite{BP}) a function $p\in C(J,X) $
is a \textit{polynomial} of degree $n$ if $p(s+mt)$ is a
polynomial in $m\in \mathbf{Z}_{+}$\ of degree at most  $n$ for
all $s,\,t\in J$ and of degree $n$ for some $s_0,\,t_0\in J$. Such
$n$ is the degree of $p$   ( $deg (p)$).
The space of all polynomials of degree at most $n$ is denoted $%
P^{n}(J,X) $. Moreover, $P(J,X)= \cup_{n=1}^{\infty} P^{n}(J,X)$,
 $P^n(J)=
P^{n}(J,\mathbf{C})$ and $P(J)= P(J,\mathbf{C})$.

 Basit and Pryde \cite{BP} gave other characterizations, on
semigroups, of the
 Domar polynomials  using difference operators $\Delta_h p(t)= p(t+h)-
 p(t)$. In particular,   $p: J \to \mathbf{C}$ is a Domar polynomial
 of degree  $n$ if and only if all differences $\Delta_h^{n+1} p$ of
 order $n+1$ vanish but $\Delta_{h_0}^{n} p (t_0)\not =0$ for some $h_0,\, t_0 \in J$.

 These polynomials were then used to study
 unbounded measurable functions $\phi : G \to X$ with
 finite Beurling spectra $sp_{w} (\phi) \subset \widehat{G}$ (see  \cite{BP}, [5, p.988 for the case $w=1$,
 $G=\mathbf{R}$]).
 Indeed, if $\phi \in L_w^{\infty} (G,X)= w L^{\infty} (G,X)$ for an appropriate
 weight $w$ and Haar measure, then  $sp_{w} (\phi)=\{\gamma_j: 1\le j\le
 n\}$ if and only if $\phi= \sum_{j=1}^n p_j \gamma_j$ for some
 polynomials $p_j\in P^n_{w} (G,X)=P^n (G,X)\cap L_{w}^{\infty} (G,X)$ (see \cite{BPB}, [9, Proposition 0.5, p. 22]).
  In a similar way, in Corollary 5.3 of \cite{BPB}, the eigenfunctions of certain
 $X$-valued convolution operators are characterized in terms of
 $X$-valued polynomials in $P^n_{w} (G,X)$. For further
 applications of polynomials, see \cite{BPJ}.

It is therefore natural to continue this study of polynomials on
semigroups. It is of particular interest to know when the space
$P^n (J)$  is finite dimensional.  For example using the theory of
symmetric matrices  we show that for some groups $P_{R} (G)$ is
properly contained in $P (G)$ (see Example 1.4), but if $P^1 (G)$
is finite dimensional then $P_{R} (G)= P (G)$ (see Theorem 2.6
(c)). This addresses the above question of Domar. In section 3, we
exhibit a class containing the compactly generated  abelian groups
for which
 $P^n (G)$ is finite dimensional for all $n\in \mathbf{Z}_+$.

\

{\center{1. \textbf{Extension  of polynomials}}\endcenter}

Throughout this section $J $ is a  subsemigroup of $G$ with
non-empty interior $J^{\circ}$ and $G=J-J$.  It is clear that
every homomorphism $\alpha: J\to \mathbf{R}$
 extends uniquely to a homomorphism
 $\alpha: G\to
\mathbf{R}$.  This implies that each Riss polynomial $p: J\to
\mathbf{C}$ extends uniquely to a Riss polynomial $p: G\to
\mathbf{C}$.  In Example 1.4, we show that $P_R^n  (J)$ may be
properly contained in $P^n (J)$. Nevertheless, as we prove in
Theorem 1.5, every polynomial $p: J\to X$ has a unique extension
to a polynomial $p: G\to X$. This last result may be of
independent interest.

 We recall some identities needed in the sequel. Let $\phi\in C(J,X)$, $t,s \in J$, $m\in \mathbf{Z}_+$.

(1.1)\qquad\qquad $\phi (t+ms) = \sum_{j=0}^m {m\choose
j}\Delta_s^j \phi(t)$.

\,

(1.2)\qquad\qquad $\,\Delta_s^m\phi (t)\,\,\,\,\, = \sum_{j=0}^m
(-1)^j {m\choose j} \phi(t+js)$.

\,

Moreover, if $\phi \in P^m (G,X)$ and $t,s \in G$ then

\,

(1.3)\qquad \qquad $\,\phi (t-s) \,\,\,\,= \sum_{j=0}^m (-1)^j
\Delta_s^j \phi(t)$.

\,

Identities (1.1), (1.2) can be verified by induction and for (1.3)
 set $z=t-s$. Then
$\Delta_s^{m+1} \phi(z)=0$ and so from (1.2),
 $\sum_{j=0}^{m+1} (-1)^j {m+1\choose j}
\phi(z+js)=0$. Thus by (1.1)
 $\phi (t-s)$$=  \sum_{j=1}^{m+1}
(-1)^{j+1} {m+1\choose j} \phi(t+(j-1)s)$
 $\,\,\,\,\,=\sum_{j=1}^{m+1} (-1)^{j+1} {m+1\choose j}\sum_{k=0}^{j-1}
{j-1\choose k}\Delta_s^k\phi(t)$

 \noindent $=  \sum_{k=0}^{m} (\sum_{i=k}^{m} (-1)^i {m+1\choose
i+1}{i\choose k}) \Delta_s^k \phi (t)=\sum_{j=0}^m (-1)^j
\Delta_s^j \phi(t)$, because $ \sum_{i=k}^{m} (-1)^i {m+1\choose
i+1}{i\choose k}$$= \sum_{i=0}^{m} (-1)^i {m+1\choose
i+1}\frac{1}{k!} \,D^k t^i|\, _{t=1} $ $\,\,\,=\frac{1}{k!} \,D^k
(t^{-1} \sum_{j=1}^{m+1} (-1)^{j+1}{m+1\choose j}\, t^j)|\,
_{t=1}$

$\qquad\qquad\qquad\qquad\qquad\qquad=\frac{1}{k!} D^k
 (t^{-1}(1-(1-t)^{m+1}))|\, _{t=1}= (-1)^k$.

\

\noindent\textbf{Examples 1.1. }

 (a) If $G$ is $\mathbf{Z}$ or $\mathbf{R}$ then $%
P^{n}(J,X) $ is the space of\ ordinary polynomials.  Indeed, that
each (ordinary) polynomial is in $P^{n}(G,X)$ is clear.
Conversely, if $\ p\in P^{n}(\mathbf{Z},X)$ then $\Delta
_{1}^{n}p(m)=c$, a constant, and $\Delta
_{t}^{n}(p(m)-cm^{n}/n!)=0$ for all $t\in \mathbf{Z.}$ An
induction argument shows $\ p$ is an (ordinary) polynomial. If $\
p\in P^{n}(\mathbf{R},X)$
then $p|_{\mathbf{Z}}\in P^{n}(\mathbf{Z},X)$ and so $p|_{\mathbf{Z}}=q|_{%
\mathbf{Z}}$ for some (ordinary) polynomial $q:\mathbf{R\rightarrow }X.$ If $%
t=a/b$ where $a,b$ are non-zero integers, then $p(t+m/b)=q(t+m/b)$
for $m\in -a+b\mathbf{Z.}$ These are polynomials in $m\in
\mathbf{Z}_{+}$ and so agree for all $m\in \mathbf{Z}_{+}$. In
particular $p(t)=q(t)$ for all rationals $ t $ and, by continuity,
for all reals $t.$

(b) If $\alpha:G\rightarrow \mathbf{C}$ is a  homomorphism, then
$\alpha\in P^{1}(G)$. Conversely, if $p\in P^{1}(G)$, then
$\alpha= p-p(0) $ is a homomorphism.

(c) If $p\in P_{w}(G,X)$ and $f\in L_{w}^{1}(G),$ then $p\ast f\in
P_{w}(G,X).$ Indeed $\Delta _{t}^{n+1}(p\ast f)=(\Delta _{t}^{n+1}p)\ast f$
for all $t\in G.$

\,

Occasionally it will be necessary for us to to assume that $
P^{n}(J)$ is finite dimensional. That this condition is not always
satisfied is shown by the following examples.

\,

\noindent\textbf{Example 1.2. } (a) Let
$G=\{s:\mathbf{N\rightarrow Z}$ ; $\ s$ has finite support$\}$
with the discrete topology. So $G$ is countable, locally compact,
$\sigma $-compact and not finitely generated. Moreover the
evaluation maps $p_{n}:G\rightarrow \mathbf{C}$ defined by
$p_{n}(s)=s(n)$ for $n\in \mathbf{ N}$ are polynomials of degree
$1$ and so $\dim \,P^{1}(G)=\infty $. Note also that $J= \{s\in G:
s(n) =0 $ for $ n < 0\} $ is a subsemigroup with $G=J-J$.

(b) The space $l^2$ of square summable real sequences with its
norm topology is not locally compact. However, if $G= l^2_d$
denotes $l^2$ with the discrete topology, then $G$ is locally
compact. If $g=(g_k)$ is a bounded real sequence and $s=(s_k)\in
l^2$ then
\smallskip

(1.4)\qquad \qquad $p(s)= \sum_{k=1}^{\infty} g_k s_k^2$, \,

\noindent is a   polynomial of degree $2$  on both $G$ and $l^2$.
Thus $P^2(G)$ and $P^2(l^2)$ are both infinite dimensional.

\,

\noindent\textbf{Lemma 1.3. } Every real valued  Riss polynomial
which is homogeneous  of order $2$ can be written in the form
$p=\sum_{j=1}^m \varepsilon_j \alpha_j^2 $ where $\alpha_j : G\to
\mathbf{R}$ are linearly independent  homomorphisms and
$\varepsilon_j=\pm 1$.

\noindent \begin{proof} A scaling argument shows that
$p=\sum_{i,j=1} ^n c_{i,j}p_i p_j$ where $p_j$ are linearly
independent  homomorphisms and $c_{i,j}$ are reals. We may also
choose $c_{i,j}=c_{j,i}$. By Sylvester's law of inertia (see [7,
p. 223]), matrix $(c_{i,j})$ is congruent to a diagonal matrix
diag$(\varepsilon_1, \cdots, \varepsilon_n)$ where
$\varepsilon_j=\pm 1$ for $1\le j\le m$ and $\varepsilon_j=0$ for
$m+1\le j\le n$, via a non-singular matrix. Hence $p=\sum_{j=1}^m
\varepsilon_j \alpha_j^2 $ where the $\alpha_j$ are  linearly
independent linear combinations of the $p_i$.
\end{proof}

\noindent\textbf{Example 1.4. }  Assume $G=l_d^2$, $g =(g_k)$ is a
bounded real  sequence, $g_k
 >0$ for all $k$ and $p$ is defined by (1.4). Then $p$ is not a Riss
 polynomial.

\noindent \begin{proof} Suppose on the contrary that $p$ is a Riss
 polynomial. Then $p=\sum_{j=1}^m \varepsilon_j \alpha_j^2 $ as in
 Lemma 1.3. The $\alpha_j$ are homomorphisms and therefore
 $\mathbf{Q}$-linear. We show that they are $\mathbf{R}$-linear.
 Note first that $p: l^2 \to \mathbf{R}$ is continuous. Moreover,
 $\Delta_h p(t)= p(h)+ \sum_{j=1}^m  2 \varepsilon_j \alpha_j (h) \alpha_j
 (t)$. Let $a=(\alpha_1, \cdots, \alpha_m): l^2\to \mathbf{R}^m$.
 Since the $\alpha_j$ are linearly independent we can find $h_i \in
 l^2$, $1\le i \le m$ such that $a(h_1), \cdots, a(h_m)$ are linearly independent in
 $\mathbf{R}^m$. Hence the linear system $\sum_{j=1}^m  2 \varepsilon_j \alpha_j (h_i) \alpha_j
 (t)= \Delta_{h_i} p(t)- p(h_i)$ can be solved  uniquely to express the $\alpha_j $ as linear combinations of
 the $ \Delta_{h_i} p- p(h_i)$. Hence each $\alpha_j: l^2 \to
 \mathbf{R}$ is continuous and since it is  $\mathbf{Q}$-linear it is $\mathbf{R}$-linear.
 Hence $a: l^2 \to \mathbf{R}^m$ has infinite dimensional kernel.
 This contradicts the fact that $p$ has only one zero.
\end{proof}

\noindent\textbf{Theorem 1.5.} Assume that $J$ is a subsemigroup
of $G$ with non-empty interior such that $G=J-J.$ Then  for each
$n$ the restriction map $ r:P^{n}(G,X)\rightarrow P^{n}(J,X)$ is a
linear bijection. In particular, every polynomial on $J$ has a
unique extension to $G$.

\noindent\begin{proof}  Firstly, let $p\in $ $P^{n}(G,X)$ be zero
on $J$. For any $t\in G$ there are $u,v\in J$ with $t=u-v.$ Now $
p(t+mv)=p(u+(m-1)v)$ is a polynomial in $m\in \mathbf{Z}_{+}$
which is zero for all $m\geq 1.$ It is therefore zero for $m=0,$
showing $p(t)=0$ and $r$ is one-to-one. Secondly, let $q\in $
$P^{n}(J,X)$. We define an extension $p$ of $q$ to $G$ as follows.
If $t\in G$ then $t=u-v$ for some $u,v\in J$ and by (1.3) it is
natural to
 set $p(t)=\sum_{j=0}^{n}(-1)^{j}\Delta _{v}^{j}q(u)$. If also $t=
\widetilde{u}-\widetilde{v}$\ where
$\widetilde{u},\widetilde{v}\in J$ we must show

\,

(1.5) $\qquad \qquad  \sum_{j=0}^{n}(-1)^{j}\Delta
_{v}^{j}q(u)=\sum_{j=0}^{n}(-1)^{j}\Delta
_{\widetilde{v}}^{j}q(\widetilde{u} ).$

\,

\noindent To do this define a function $L:P^{n}(J,X)\rightarrow X$
by $ L(q)=q(u)-\Delta _{v}q(\widetilde{u})=q(\widetilde{u})-\Delta
_{\widetilde{v} }q(u).$ We prove by induction on $k\in
\mathbf{Z}_{+}$ that

\,

(1.6)  $ \qquad\sum_{j=0}^{n}(-1)^{j}\Delta
_{v}^{j}q(u)=\sum_{j=0}^{k-1}L(\Delta _{v}^{j}\Delta
_{\widetilde{v}}^{j}q)+\Delta _{v}^{k}\Delta _{\widetilde{v}
}^{k}q(w)$ \,\,  if \, $q$ is of degree $n=2k,$

 \,

 \qquad\qquad $\sum_{j=0}^{n}(-1)^{j}\Delta
_{v}^{j}q(u)=\sum_{j=0}^{k}L(\Delta _{v}^{j}\Delta
_{\widetilde{v}}^{j}q)$ \,\qquad  if \,\, $q$ is of degree
$n=2k+1$

\,

\noindent where $w$ is an arbitrary element of $J.$ When $n=0,$
$q(u)=q(w)$ and when $ n=1$, $\Delta _{v}q(u)=\Delta _{v}q(w)$ and
so $q(u)-\Delta _{v}q(u)=q(u)-\Delta _{v}q(\widetilde{u})=L(q).$
Hence, (1.6) holds for $k=0.$ Assume (1.6) holds for polynomials
of degree less than $2k$. If $n=2k,$ we can apply (1.6) to $\Delta
_{v}q$ and obtain

\qquad $\sum_{j=0}^{n}(-1)^{j}\Delta
_{v}^{j}q(u)=q(u)-\sum_{j=0}^{k-1}L(\Delta _{v}^{j}\Delta _{\widetilde{v}%
}^{j}\Delta _{v}q)$

 \qquad $
=q(u)+\sum_{j=0}^{k-1}\Delta _{\widetilde{v}}\Delta _{v}^{j}\Delta _{%
\widetilde{v}}^{j}\Delta _{v}q(u)$
$ -\sum_{j=0}^{k-1}\Delta _{v}^{j}\Delta _{%
\widetilde{v}}^{j}\Delta _{v}q(\widetilde{u})$

\qquad $%
=\sum_{j=0}^{k}\Delta _{v}^{j}\Delta _{\widetilde{v}}^{j}q(u)-%
\sum_{j=0}^{k-1}\Delta _{v}^{j}\Delta _{\widetilde{v}}^{j}\Delta _{v}q(%
\widetilde{u})$
 $%
=\sum_{j=0}^{k-1}L(\Delta _{v}^{j}\Delta
_{\widetilde{v}}^{j}q)+\Delta _{v}^{k}\Delta
_{\widetilde{v}}^{k}q(u).$

If $n=2k+1,$ we can apply this last result to $\Delta _{v}q$ and
obtain

\qquad $\sum_{j=0}^{n}(-1)^{j}\Delta
_{v}^{j}q(u)=q(u)-\sum_{j=0}^{k-1}L(\Delta _{v}^{j}\Delta _{\widetilde{v}%
}^{j}\Delta _{v}q)-\Delta _{v}^{k}\Delta
_{\widetilde{v}}^{k}\Delta _{v}q(w)$

\qquad $%
=q(u)+\sum_{j=0}^{k-1}\Delta _{\widetilde{v}}\Delta _{v}^{j}\Delta _{%
\widetilde{v}}^{j}\Delta _{v}q(u)-\sum_{j=0}^{k-1}\Delta _{v}^{j}\Delta _{%
\widetilde{v}}^{j}\Delta _{v}q(\widetilde{u})-\Delta _{v}^{k}\Delta _{%
\widetilde{v}}^{k}\Delta _{v}q(\widetilde{u})$

\qquad $%
=\sum_{j=0}^{k}\Delta _{v}^{j}\Delta _{\widetilde{v}}^{j}q(u)-\sum_{j=0}^{k}%
\Delta _{v}^{j}\Delta _{\widetilde{v}}^{j}\Delta
_{v}q(\widetilde{u}) =\sum_{j=0}^{k}L(\Delta _{v}^{j}\Delta
_{\widetilde{v}}^{j}q).$

\noindent Hence, (1.6) is proved and (1.5) follows, which means
$p$ is well-defined.
 Note that if $t\in J$ then we can take $t=u$, $v=0$ and so $p(t)=q(t)$. So $%
p $ is an extension of $q$. Next, we show $p$ is continuous. Since
$J$ has an interior point $s_{0}$, there is an open neighborhood
$W$ of $0$ in $G$
such that $s_{0}+W\subset J$.  Moreover, if $t=u-v$, where $u,v\in J$ then $t=%
\widetilde{u}-\widetilde{v}$, where $\widetilde{u}=s_{0}+u,\widetilde{v}%
=s_{0}+v$. Let $(t_{\alpha })$ be a net in $G$ converging to $t$.
We may
suppose $t_{\alpha }=t+w_{\alpha }$ where $w_{\alpha }\in W$. Setting $%
u_{\alpha }=\widetilde{u}+w_{\alpha }$ and $v_{\alpha
}=\widetilde{v}$ we
find $u_{\alpha }, v_{\alpha }\in J, t_{\alpha }= u_{\alpha }-v_{\alpha }$ and $%
(u_{\alpha })\rightarrow \widetilde{u}$. So $p(t_{\alpha
})\rightarrow p(t)$.
Finally, if $t_{j}=u_{j}-v_{j}$ where $u_{j},v_{j}\in J$ and $m\in \mathbf{Z}%
_{+}$ then $p(t_{1}+mt_{2})=\sum_{j=0}^{n}(-1)^{j}\Delta
_{v_{1}+mv_{2}}^{j}q(u_{1}+mu_{2})$ which is a polynomial in $m$
of degree at most $n$. So $p\in $ $P^{n}(G,X),$ proving $r$ is
onto. \end{proof}

{\center{\textbf{2. Sufficient conditions for all polynomials to
be Riss}
\endcenter}}

In this section we develop a method enabling us to determine the
dimension of the space of polynomials $P^n (J)$ and  to  give
conditions under which $P^n (J) =P_R^n (J)$.

\noindent\textbf{Lemma 2.1. } Assume that $\phi\in P^n (J,X)$ so
$\phi (t+ms)= \sum_{j=0}^n  a_j (t,s)m^j$ for all $t, s\in J$,
$m\in \mathbf{Z}_+$. Then $a_n(t,s)= a_n (0,s)$.

\noindent\begin{proof}  Applying (1.2) to the polynomial on
$\mathbf{R}$ given by $p_k (u)=u^k$ where $0\le k \le n$ and using
$\Delta_1^n p_k (0)=  n!\delta_{k,n}$ we get

(2.1)\qquad $ \sum_{j=0}^n (-1)^j {n\choose j} j^k =
n!\delta_{k,n}$ for  $0\le k \le n$ ($0^0=1$).

\noindent Now consider $\phi \in P^n (J,X)$. We have
$\Delta_s^n\phi (t)$ is independent of $t$. By (1.2) again

$\Delta_s^n\phi (t)= \sum_{j=0}^n (-1)^j {n\choose j} \phi(t+j s)=
\sum_{j=0}^n (-1)^j {n\choose j} \sum_{k=0}^n a_k(t,s)j^k$

$=\sum_{k=0}^n a_k(t,s)\sum_{j=0}^n (-1)^j {n\choose j} j^k
=\sum_{k=0}^n a_k(t,s)n!\delta_{k,n} \,\,= n! a_n (t,s)$.

\noindent So $ a_n (t,s)= a_n (t,0)$.
\end{proof}

\noindent\textbf{Lemma 2.2. } If  $\phi\in P^n (J,X)$ and $\phi
(ms)= m^k \phi (s)$ for some $0 \le k < n$ and all  $m\in
\mathbf{Z}_+$, $s\in J$ then $\phi\in P^k (J,X)$.

\noindent\begin{proof} From (1.1), we have $\phi (t+ms)=
\sum_{j=0}^n a_j(t,s)\, m^j$,  so $\phi (ms)= \sum_{j=0}^n
a_j(0,s)\, m^j$ $= m^k\, \sum_{j=0}^n a_j(0,s)$.
 It follows that $a_j (0,s)=0$ for all $j\not =k$. By Lemma 2.1,
$a_n (t,s)= a_n(0,s)=0$. It follows $\phi\in P^{n-1} (J,X)$.
Applying this argument repeatedly we obtain  $\phi\in P^k (J,X)$.
\end{proof}

\noindent\textbf{Proposition 2.3. } If  $\phi\in P^n (J,X)$ then
$\phi=\sum_{j=0}^n a_j$ where $a_j \in P^j (J,X)$ and $a_j(mt)=
m^j a_j (t)$ for all $m\in \mathbf{Z}_+$, $t\in J$.

\noindent\begin{proof} From (1.1), we have $\phi (mt) =
\sum_{j=0}^m {m\choose j}\Delta_t^j \phi(0)$ for all $m\in
\mathbf{Z}_+$. But $\Delta_t^j \phi=0$ if $j\ge n+1$ and so $\phi
(mt) = \sum_{j=0}^n {m\choose j}\Delta_t^j \phi(0)= \sum_{j=0}^n
a_j(t)\, m^j$, where $a_j (t)= \sum_{k=0}^n a_{jk} \Delta_t ^k
\phi(0)$ for some $a_{jk}\in \mathbf {R}$. Thus $a_j\in P^n
(J,X)$. Moreover, $\phi (mlt)=\sum_{j=0}^n a_j (lt)m^j$. So, $a_j
(lt)= l^j a_j (t)$ for all $l\in \mathbf {Z}_+$. By Lemma 2.2,
$a_j \in P^j (J,X)$.
\end{proof}

  Recall that a group is called a $torsion\,\, group$
 if every element has finite order. If the orders of the elements
 are bounded the group is said to be of $bounded\,\, order$.
 A group  is
$torsion\,\, free$ if no element other than  the identity is of
finite order (see [6,  Definition, p. 88]).

\,

 \noindent\textbf{Lemma 2.4}. Let $K$ be a  compact or torsion abelian topological group.
Then $P(K,X)= P^0(K,X)$.

\noindent\begin{proof}\, Let $p\in P(K,X)$.  If   $K$ is compact,
$p (K)$ is compact and if $K$ is torsion, then the sequence  $ (p
(t+ms))$, $m\in \mathbf{Z}_+$ is compact for each $t,s \in K$. If
$p$ is not constant then for some $t,s \in K$ one has $p(t+ms)$ is
an unbounded polynomial in $m\in \mathbf{Z}_+$. This is a
contradiction which proves that $deg \,(p)=0$.
\end{proof}

 We give the proof of the following
needed result stated without a proof in \cite{MWG}.

\,

 \noindent\textbf{Theorem 2.5 (Mackey)}. The following
statements are equivalent.

(a) \, For each $t\in G\setminus \{0\}$ there is a homomorphism
$\alpha: G\to \mathbf{R}$ such that $\alpha(t)\not = 0$.

(b)\, $G= \mathbf{R}^m \times F$ for some $m\in \mathbf{Z}_+$  and
some discrete torsion free group $F$.

\noindent\begin{proof} $(a)\Rightarrow (b)$: By the  principal
structure theorem [12,Theorem  2.4.1], $\widehat{G}$ has an open
and closed subgroup  $\Gamma_0 = \mathbf{R}^m \times K$, where $K$
is compact. The annihilator  $H= \Gamma_0^{\perp}$ is isomorphic
to the dual of the discrete group $\widehat{G}/\Gamma_0$ (see [12,
Theorem 2.1.2, p. 35]). So $H$ is compact and thus every
homomorphism $\alpha: G\to \mathbf{R}$ is zero on $H$ (see also
Lemma 2.4). By (a) we conclude $H= \{0\}$ and therefore
$\widehat{G}=\Gamma_0 = \mathbf{R}^m \times K$. Hence $G=
\mathbf{R}^m \times F$ where $F= \widehat{K}$ is discrete. By (a)
we also conclude that $G$ and therefore $F$ are torsion free.

$(b)\Rightarrow (a)$: If $t\in F\setminus \{0\}$ then using Zorn's
lemma we may establish the existence of a homomorphism $\alpha:
F\to \mathbf{R}$ with $\alpha (t)\not =0$ (see also [12, Theorem
(Kaplansky), p. 44]). Then $\alpha$ extends readily to a
homomorphism $\alpha: G\to \mathbf{R}$. Now (a) follows.
\end{proof}

In the sequel for $Y=X$ or $\mathbf{R}$ we use

 \,

\qquad  $H_0 (Y)= \{x\in G: \alpha (t)= 0$ for all homomorphisms
$\alpha:G\to Y\}$ and $H_0= H_0 (\mathbf{R})$.

\,

Clearly $H_0 (X)\subset H_0$. We claim that if $X\not = \{0\}$,
then $H_0 (X)=H_0$.
 Indeed this is clear if $X$ is locally convex, by the Hahn-Banach
theorem. Assume $X$ is not locally convex and there is $ t\in H_0
$ such that $t\not \in
 H_0 (X)$. Then there is a homomorphism
$\alpha:G\to X$ such that $y=\alpha (t)\not = 0$. Let $X_y =
\mathbf{C}\, y$ and $P_y: X\to X_y $ be a projection map. Then
$\beta
 =P_y\circ\alpha: G\to X_y  $ is a homomorphism with $\beta (t)=
\alpha (t) \not = 0$. But $X_y$ is locally convex. This implies
$\alpha (t)=0$, a contradiction which proves our claim.

\,

\noindent\textbf{Proposition 2.6}. If $\phi\in P^n (G,X)$, then
$\phi(t)=\phi(0)$ for all $t\in H_0$.

\noindent\begin{proof} If $n=1$, $\phi(t)-\phi(0)=\alpha (t)$,
$\alpha: G\to X$ a homomorphism. So $\phi(t)=\phi(0)$ on $H_0$ for
all $\phi\in P^1 (G,X)$. As induction hypothesis  assume $\phi(x)=
\phi (0)$ on $H$ for all $\phi\in P^{n-1} (G,X)$. By Proposition
2.3, it suffices to prove the claim for all $\phi\in P^n (G,X)$
with $\phi (mt)= m^n \phi (t)$. But $\Delta_s\phi \in P^{n-1}
(G,X) $ so $\Delta_s\phi (t)=\Delta_s\phi (0) $ for $t\in H_0$. In
particular $\Delta_{mt}\phi (t) =\Delta_{mt}\phi (0) $ for $t\in
H_0$, $m\in \mathbf{Z}_+$. Thus $ \phi ((m+1)t)- \phi (t)=\phi (m
t)- \phi (0)$, that is $(m+1)^n \phi (t)- \phi (t)= m^n \phi (t)
$. Hence $((m+1)^n - (m^n +1))\phi (t)=0$. But $n >1$, so $\phi
(t)=0$ as required.
\end{proof}

\noindent\textbf{Theorem 2.7}. With $H_0$ as above:

 (a)\, $ P^n (G/H_0,X)=P^n (G,X)$ for $n\in \mathbf{Z}_+$.

(b)\, $G/H_0= \mathbf{R}^m\times F$ for some $m\in \mathbf{Z}_+$
and some torsion free discrete group $F$.

(c)\,  If $P^1(G)$ is finite dimensional, then  $G/H_0=
\mathbf{R}^m\times \mathbf{Z}^k$ for some  $ k, m\in \mathbf{Z}_+$
and $P ^n (G)= P_R^n (G)$ for all $ n\in \mathbf{Z}_+$.

\noindent\begin{proof} (a)\, Let $\pi: G\to G/H_0$ be the
canonical projection. By Proposition 2.6 the map $ P^n
(G/H_0,X)\to P^n (G,X)$ defined by $\phi\to \phi\circ \pi$ is an
isomorphism.

(b)\, Let  $\xi =\pi (t) \in G/H_0$, $t\not \in H_0$.
 By definition
of $H_0$ there is a homomorphism $\alpha :G\to \mathbf{R}$ with
$\alpha(t) \not = 0$. Its preimage under the mapping in (a) with
$X=\mathbf{C}$ is a homomorphism $a: G/H_0\to \mathbf{R} $
 with $a(\xi)\not =0$.  Thus
 by  Theorem 2.5,  $G/H_0= \mathbf{R}^m\times
F$ for some $m \in \mathbf{Z}_+$ and some torsion free discrete
group $F$.

(c)  Let $F_0 \subset F$ be any finitely generated subgroup of
$F$. By [6, Theorem 9.3, p. 90] $F_0= \mathbf{Z}^k$ for some $k\in
\mathbf{Z}_+$. The coordinate functions $\alpha_j: F_0\to
\mathbf{R}$ extend  by [12, Theorem (Kaplansky), p. 44] to
homomorphisms $\alpha_j: F\to \mathbf{R}$. By part (b) they extend
further to linearly independent homomorphisms $\alpha_j: G/H_0\to
\mathbf{R}$. Since $P^1(G)$ is finite dimensional we conclude that
there is a maximal such $k\in \mathbf{Z}_+$. Hence $F=
\mathbf{Z}^k$ and  $ P ^n(G)=P ^n(G/H_0)= P_R^n (G/H_0)=P_R^n
(G)$.
\end{proof}

\noindent\textbf{Corollary 2.8}. If $G=J-J$ where $J^{\circ}\not =
\emptyset$ and $P^1 (J)$ is finite dimensional, then $G/H_0=
\mathbf{R}^m\times \mathbf{Z}^k$ for some $m, k\in \mathbf{Z}_+$
and  $P^n (J)$ is finite dimensional for each $n\in \mathbf{Z}_+$.

\begin{proof} By Theorem 1.5, $P^1 (G)$ is finite
dimensional. Therefore the  result follows from Theorem 2.7 (c).
\end{proof}

{\center{\bf{3. A class of groups with $P^n (G)$ finite
dimensional}}\endcenter}

In this section we  assume  also that $X$ is locally convex so
that the Hahn-Banach theorem may be applied. We exhibit a general
class of groups  $G$ under which $P^n (G)$ is finite dimensional.

\,

\noindent \textbf{Lemma 3.1.}  (a)\, If $P^{n}(J)$ is finite
dimensional
then $P^{n}(J,X)=P^{n}(J)$ $\otimes X$.

(b) If $%
P_{w}^{n}(G)$ is finite dimensional then $%
P_{w}^{n}(G,X)=P_{w}^{n}(G)$ $\otimes X$.

\noindent\begin{proof} (a) \, Clearly, $P^{n}(J)$ $\otimes X\subseteq $ $%
P^{n}(J,X).$ For the converse, which is clearly true when $n=0,$\
we use induction on $n.$\ Let $\{p_{1},...,p_{k}\}$ be a basis of
a complement $Q$\ of $P^{n-1}(J)$ in$\ P^{n}(J).$ Since $\deg
(p_{1})=n$ we can choose $t_{1}\in J$ such that $\Delta
_{t_{1}}^{n}p_{1}\neq 0$. Since each $\Delta _{t_{1}}^{n}p_{j}$ is
a constant we can set $q_{1}=p_{1}/\Delta
_{t_{1}}^{n}p_{1}$ and choose $\lambda _{1}\in \mathbf{C}$\ such that $%
\Delta _{t_{1}}^{n}(p_{2}-\lambda _{1}q_{1})=0$. Set
$q_{2}=p_{2}-\lambda
_{1}q_{1}$. Then $\deg (q_{2})=n$ for otherwise $q_{2}\in Q\cap P^{n-1}(J,%
\mathbf{C})=\{0\},$ contradicting the linear independence of
$p_{1},p_{2}$. Hence we can choose $t_{2}\in J$ such that $\Delta
_{t_{2}}^{n}q_{2}\neq 0$. Continuing in this way, we obtain a
basis $\{q_{1},...,q_{k}\}$ of $Q$ and a subset
$\{t_{1},...,t_{k}\}$ of $J$ such that $\Delta
_{t_{i}}^{n}q_{j}=\delta _{i,j}$. Now let $p\in $ $P^{n}(J,X)$. For each $%
x^{\ast }\in X^{\ast }$\ we have $x^{\ast }\circ p\in P^{n}(J)$.
Hence $x^{\ast }\circ p=\sum_{j=1}^{k}q_{j}c_{j}(x^{\ast })+r (x^*)$ for some $%
c_{j} (x^*)\in \mathbf{C}$ and $r (x^*)\in P^{n-1}(J)$. But
$x^{\ast }\circ \Delta _{t_{i}}^{n}\,p=\Delta _{t_{i}}^{n}(x^{\ast
}\circ p)=c_{i}(x^{\ast })$ and so $c_i (x^*)=  x^* (c_i)$ where
$c_i = \Delta^n _{t_{i}} \, p$. Moreover, $x^{\ast }\circ
(p-\sum_{j=1}^{k}q_{j}c_{j})=r (x^*)\in P^{n-1}(J)$ and so by the
Hahn-Banach theorem
$\Delta^n_t (p-\sum_{j=1}^{k}q_{j}c_{j})=0$  for all $t\in J $  showing $%
p-\sum_{j=1}^{k}q_{j}c_{j}\in P^{n-1}(J,X)$. By the induction hypothesis $%
P^{n-1}(J,X)=P^{n-1}(J)$ $\otimes X$ and hence $p\in P^{n}(J,%
\mathbf{C})$ $\otimes X$ as required.

(b) This is proved in the same way as above.
\end{proof}

\noindent\textbf{Lemma 3.2. }\, Let $J_1,\, J_2$ be abelian
topological semigroups. If $P^{n}(J_{1})$ is finite dimensional,
then

 $ \qquad \qquad \qquad P^{n}(J_{1}\times J_{2})=\,
 \sum_{m=0}^{n}P^{m}(J_{1})\otimes
P^{n-m}(J_{2})$.

\noindent\begin{proof}\ The inclusion $\supseteq $ is clear. For
the converse, take any $p\in P^{n}(J_{1}\times J_{2})$. As in the
proof of
Lemma 3.1, we can find for each $m=1,...,n$ a basis $%
\{q_{1}^{m},...,q_{k_{m}}^{m}\}$ of a complement of $P^{m-1}(J_{1}%
)$ in$\ P^{m}(J_{1})$ and a subset $\{s_{1}^{m},...,s_{k_{m}}^{m}%
\}$ of $J_{1}$ such that $\ \Delta
_{s_{i}^{m}}^{m}q_{j}^{m}=\delta _{i,j}$.
Also let $\left\{ q_{1}^{0}\right\} $ be a basis of $P^{0}(J_{1}%
). $ For any $t\in J_{2}$ we have $p(.,t)\in P^{n}(J_{1})$ and so
$p(s,t)=\sum_{m=1}^{n}\sum_{j=1}^{k_{m}}q_{j}^{m}(s)r_{j}^{m}(t)$ for some $%
r_{j}^{m}(t)\in \mathbf{C}$. We prove by backward induction on
$h$\ that $r_{j}^{h}\in P^{n-h}(J_{2})$. Now

$\Delta
_{(s_{i}^{n},0)}^{n}p(s,t)=\sum_{m=1}^{n}\sum_{j=1}^{k_{m}}\Delta
_{s_{i}^n}^n q_{j}^{m}(s) r_{j}^{m}(t)=r_{i}^{n}(t)$

\noindent so each $r_{i}^{n}$ is a constant as required. So
suppose each $ r_{j}^{m}\in P^{n-m}(J_{2})$ for
 $n\geq m\geq h+1\,\,\,$ and
$\,\,\, 1\leq
j\leq k_{m}$. Then

$p(s,t)-\sum_{m=h+1}^{n}%
\sum_{j=1}^{k_{m}}q_{j}^{m}(s)r_{j}^{m}(t)=\sum_{m=0}^{h}\,
\sum_{j=1}^{k_{m}}q_{j}^{m}(s)r_{j}^{m}(t)\, \in \, P^{n}(J_{1}\times J_{2},%
\mathbf{C})$ and

 $\Delta
_{(s_{i}^{h},0)}^{h}\sum_{m=0}^{h}%
\sum_{j=1}^{k_{m}}q_{j}^{m}(s)r_{j}^{m}(t)=r_{i}^{h}(t)$. So each $%
r_{i}^{h}\in P^{n-h}(J_{2})$

\noindent and the proposition is proved. \end{proof}

 \noindent\textbf{Lemma 3.3}.  Let $H$ be a
 closed subgroup of $G$ such that $G/H$ is a torsion group.

 (a) The restriction map $r:$ $P^{n}(G,X)\rightarrow P^{n}(H,X)$ is
 one-to-one and so $\dim $ ($P^{n}(G,X))\leq $ $\dim (P^{n}(H,X))$.

 (b) If also $G/H$ is of bounded order and $P^{n}(H)$ is finite
 dimensional, then

 $r:$ $P^{n}(G,X)\rightarrow P^{n}(H,X)$ is a linear isomorphism.

 \noindent\begin{proof}\ (a) Let $p\in P^{n}(G,X)$ satisfy $p(t)=0$ for all $t\in H$.
 \ Let $s\not \in  H$ and let $\pi :G\rightarrow G/H$ be the
 quotient map. Since $G/H$ is a torsion group,  $\pi (ks)=0$, meaning  $ks\in H$, for some $k\in
 \mathbf{N.}$ Hence $p(mks)=0$
 for all $m\in \mathbf{N.}$ But $p(ms)$ is a polynomial in $m\in \mathbf{Z}
 _{+}$ and so is zero. In particular, $p(s)=0$ showing $r$ is
 one-to-one.

 (b) Let $\{p_{1},...,p_{m}\}$ be a basis of $P^{n}(H)$
 and choose $k\in \mathbf{N}$ such that\ $kt\in H$ for all $t\in G$. Define $
 q_{j}(t)=p_{j}(kt)$ and suppose $\sum_{j=1}^{m}c _{j}q_{j}=0$ on
 $G$ for some $c _{j}\in \mathbf{C.}$ Then $\sum_{j=1}^{m}c
 _{j}p_{j}=0$ on $kG$. But $kG$ is a closed subgroup of $H$ such
 that $H/kG$ is a torsion
 group. By part (a), $\sum_{j=1}^{m}c _{j}p_{j}=0$ on $H$. Hence each $
 c _{j}=0,$ showing $\{q_{1},...,q_{m}\}$ is linearly independent and $
 \dim $ ($P^{n}(G))\geq $ $\dim (P^{n}(H))$. Therefore $
 r:$ $P^{n}(G,X)\rightarrow P^{n}(H,X)$ is a linear isomorphism, by (a) when $
 X=\mathbf{C}$ and then by Lemma 3.1 for general $X$.
 \end{proof}

Reiter, [9, p.142], introduced the notion of $w$-$uniform\,
continuity$ for functions $\phi\in C(G,X)$, namely sup $_{t\in G}
||\Delta_h \phi (t)/ w(t)|| \to 0$ as $h\to 0$ in $G$. The space
of all such functions for  which   $\phi/w$  is bounded is denoted
by $BUC_w (G,X)$. The weight $w$ has $polynomial\,\, growth$ of
order $N\in \mathbf{Z}_+$ if it satisfies conditions (3.1), (3.2)
of \cite{BP}.

\,

\noindent\textbf{Theorem 3.4}. Assume $G$  has an open subgroup $G_0$ such that $%
G_0=\mathbf{R}^{m}\times K$ for some $m\in \mathbf{Z}_+$ and some
compact group $K$, and $G/G_0=\mathbf{Z}^{k}\times F$ for some
$k\in \mathbf{Z}_+$ and some torsion group $F$. Then for each
$n\in \mathbf{Z}_+$,

 (a) $P^{n}(G)$ is finite dimensional.

\noindent Assume also $F$ is of bounded order.

(b)  For each $p\in P^{n}(G,X)$ there exist $p_{j}\in P^{n}(G,X)$ and $%
q_{j}\in P^{n}(G)$\ with $q_{j}(0)=0$ such that $\Delta
_{h}p(t)=\sum_{j=1}^{k}p_{j}(t)q_{j}(h)$\ for all $h,t\in G$.

(c) $P_{w}(G,X)\subseteq BUC_{w}(G,X)$.

(d) If also $w$ has polynomial growth $N$, then $P_w(G)$ is finite
dimensional.

\noindent\begin{proof} (a) Note that $G_0$  is also closed (see
[12, Appendix B5]). Hence $G/G_0$ is discrete. Since
$G/G_0=\mathbf{Z}^{k}\times F$
 and
$\mathbf{Z}^{k}$ is a free group,  $G$ has a subgroup $H$
isomorphic to $\mathbf{Z}^{k}$ with $H \cap G_0 =\{0\}$. Since
$G_0$ is open, $G_1=G_0 \times H$ is open too.  This implies $G_1$
is also closed.  Moreover, $G/G_{1}=F$. By Lemma 3.3 (a) $\dim
P^{n}(G)\le \dim P^{n}(G_1)$. By  Examples 1.1 (a), Lemma 2.4 and
Lemma 3.2  we conclude that $\dim P^{n}(G)$ is finite.

(b) Since (b) holds when $G=\mathbf{R} ^{m}\times \mathbf{Z}^{k}$
it holds for  $G_1$. Since $F$ has bounded order, (b) follows by
Lemma 3.3 (b).

(c) If $p\in P_{w}(G,X)$ choose $p_{j},q_{j}$ as in (b).\ Since
$|| \Delta
_{h}p(t)|| \le cw(t)\sum_{j=1}^{k} || q_{j}(h)|| $, where $%
c=  \sup_{j} \sup _{t\in G}\,|| p_{j}(t)/w(t) ||$, it follows that
$p\in BUC_{w}(G,X)$, proving (c).

(d) The assumptions imply $P_{w}(G,X) \subset P^N(G,X)$ and
therefore the result follows from (a).
 \end{proof}

\noindent\textbf{Remarks 3.5}. (a) \, If  $G$  is compactly
generated, then $G$ satisfies all the assumptions of Theorem 3.4.
   Indeed, by the principal structure theorem for locally
compact abelian groups (see [12, Theorem 2.4.1]), each such group
$G$ has an open subgroup $G_0=\mathbf{R}^{m}\times K$ for some
$m\in \mathbf{Z}_+$ and some compact group $K$. Hence $G/G_{0}$ is
discrete. Since $G$ is compactly generated  so is $G/G_{0}$ and
therefore $G/G_{0}$ is finitely generated. So
$G/G_{0}=\mathbf{Z}^{k}\times F$ for some $k\in \mathbf{Z}_+$ and
some finite group $F$ (see Theorem 9.3 in \cite{F}).

(b) There are groups $G$ satisfying all the assumptions of Theorem
3.4 which  are not compactly generated. For example, take
$F=\mathbf{Z}_2^{\mathbf{R}} =\{ f: \mathbf{R}\to\mathbf{Z}_2 \}$,
where $\mathbf{Z}_2$ is the group of  order $2$, endowed with the
discrete topology. Then $F$ is a torsion group of bounded order
$2$ which is not compactly generated.

(c) There are locally compact groups $G$ for which $P^n (G)$ is
infinite dimensional for each $n > 0$.
 See Examples 1.2 (b).

\end{document}